\documentclass[12pt,leqno,british]{amsart}

\newif\ifpdf
\ifx\pdfoutput\undefined\pdffalse\else\pdfoutput=1\pdftrue\fi\newif\pdf
\ifpdf\relax\else
\usepackage[T1]{fontenc}
\fi

\usepackage{babel}
\usepackage{eucal}
\usepackage{url}
\usepackage{stmaryrd}
\usepackage[pagebackref]{hyperref}
\usepackage{amsfonts}
\usepackage{amsthm}
\usepackage{amsmath}
\usepackage{amscd}
\usepackage{amssymb}

\def\CC {{\mathbb C}}     
\def\RR {{\mathbb R}}     

\def\lo  {\longmapsto}

\def\mk {\mathfrak}

\def\be  {\begin{eqnarray}}
\def\ee  {\end{eqnarray}}
\def\ben {\begin{eqnarray*}}
\def\een {\end{eqnarray*}}

\def\bpr {\begin{proof}[Proof]}
\def\epr {\end{proof}}
\def\bsp {\begin{split}}
\def\esp {\end{split}}
\def\bprr {\begin{proof}[solution]}
\def\bpru {\begin{proof}[hint]}
\def\bpro {\begin{proof}[answer]}
\def\bcd {\begin{CD}}
\def\ecd {\end{CD}}

\newcommand{\scal}[1]{\left\langle#1\right\rangle}
\newcommand{\norm}[1]{\left\Vert#1\right\Vert}
\newcommand{\sco}[1]{\left(#1\right)}

\newtheorem{theorem}{Theorem}[section]

\newtheorem{prop}[theorem]{Proposition}
\newtheorem{coro}[theorem]{Corollary}
\newtheorem{zab}{Remark}[section]
\newtheorem{df}{Deffinition}[section]

\newtheorem*{ack}{Acknowledgements}

\begin{document}
\setcounter{page}{1}

\title[Complex product structures on some simple Lie groups]%
{Complex product structures on some simple Lie groups}

\author{Stefan Ivanov}
\address[Ivanov]{University of Sofia "St. Kl. Ohridski"\\
Faculty of Mathematics and Informatics,\\ Blvd. James Bourchier
5,\\ 1164 Sofia, Bulgaria,}
\email{ivanovsp@fmi.uni-sofia.bg}

\author{Vasil Tsanov}
\address[Tsanov]{University of Sofia "St. Kl. Ohridski"\\
Faculty of Mathematics and Informatics,\\ Blvd. James Bourchier
5,\\ 1164 Sofia, Bulgaria,}
\email{tsanov@fmi.uni-sofia.bg}

\begin{abstract}
We construct invariant complex product (hyperparacomplex, indefinite
quaternion) structures on the manifolds underlying the real
noncompact simple Lie groups $SL(2m-1,\RR)$, $SU(m,m-1)$ and
$SL(2m-1,\CC)^\RR$. We show that on the last two
series of groups some of these structures are compatible with
the biinvariant Killing metric. Thus we also provide
a class of examples of compact (neutral) hyperparahermitean,
non-flat Einstein manifolds.

\end{abstract}

\maketitle

\section{Introduction}

In the present paper we construct explicit examples of homogeneous
complex product structures on real semisimple Lie groups and some
related manifolds. Our examples fall into the class of CPS
known also as hyperparacomplex structures.
The relevant definitions are provided in Section \ref{prel}.

The main result is
\begin{theorem}\label{cps1}
For each $m >1$, the manifolds underlying the
Lie groups $SL(2m-1,\RR)$ and $SU(m,m-1)$ have complex
product structures, which are invariant by left translations.
\end{theorem}

The proof of the above theorem will be given in Section
\ref{stru}, there we construct explicitely a family of complex
product structures depending on parameters in each of the above Lie
groups. In this introduction we indicate some more examples of
complex product manifolds, obtained directly from Theorem
\ref{cps1}, and discuss properties of their complex product
structures. For the sake of brevity, we sometimes speak of complex
product structures on Lie algebras, rather than on the manifolds
underlying the respective Lie groups. Reduction of the latter to
the former is explained in Section \ref{prel} (see also e.g.
\cite{AnSa1}).

To our knowledge (see also Andrada and Salamon  \cite{AnSa1} p.2)
Theorem \ref{cps1} provides the first examlpes of homogeneous
complex product structures on semisimple real Lie groups.
There are many known complex product structures on solvable
groups. Andrada and Salamon \cite{AnSa1} have constructed invariant
complex product structures on the reductive groups $GL(2m,\RR)$.
From the proof of Theorem \ref{cps1}, we shall obtain complex
product structures on the groups $GL(2m,\RR)$,
embedding\footnote{see Remark \ref{rema1}.}
them as (complex product) submanifolds of $SL(2m+1,\RR)$
(i.e. inheriting the complex product structure of
Theorem \ref{cps1}). In the same way we obtain complex
product structures on $U(m,m)$.

Let $\mk{g}$ be a real Lie algebra, and let $(P,J)$ be a
complex product structure on $\mk{g}$. It is trivial to
check, that if we extend $P, J$ to the complexification
$\mk{g}^\CC$ by linearity ($P(iX) \doteq iP(X),\quad
J(iX) \doteq iJ(X)$), then we get a complex product
structure on the real Lie algebra $\sco{\mk{g}^\CC}^\RR$. Thus
we have

\begin{coro}\label{cps12}
For each $m >1$, the manifold underlying $(SL(2m-1,\CC))^\RR$ admits
homogeneous complex product structures.
\end{coro}

To get a more interesting and subtle result we combine
Theorem \ref{cps1} with Theorem 3.3 of Andrada and
Salamon \cite{AnSa1} where they prove, that a complex product
structure on a real Lie algebra $\mk{g}$ induces a hypercomplex
structure on $\mk{g}^\CC$. As $sl(2m-1,\CC) = (sl(2m-1,\RR))^\CC =
((su(m,m-1))^\CC)$ we get
\begin{coro}\label{cps7}
For each $m >1$, the manifold underlying $SL(2m-1,\CC)$ admits
homogeneous hypercomplex structures.
\end{coro}

Let $\mk{g}$ be a real Lie algebra with a complex product
structure $(P,J)$. Each nondegenerate symmetric bilinear
form $g(X,Y)$ on $\mk{g}$, satisfying the condition (\ref{me1}),
defines an invariant pseudoriemannian metric (neutral)
on the manifold $G$ (the simply connected group whose Lie
algebra is $\mk{g}$), which is compatible\footnote{See Definition
\ref{me7}.} with $(P,J)$. Thus the homogeneous complex
product structures of Theorem \ref{cps1} and Corollary
\ref{cps12} admit compatible left invariant metrics.
On $SL(2m-1,\RR)$ the Killng form is not of
neutral signature, so it cannot be compatible with any
complex product structure. On the other hand we have

\begin{theorem}\label{kill1}
For each $m>1$ the group manifolds $SU(m,m-1)$ and
$(Sl(2m-1,\CC))^\RR$ admit complex product structures, which are
compatible with the biinvariant metric induced
by the respective Killing forms.
\end{theorem}
The proof will also be given in Section \ref{stru}.

Let $(P, J)$ be a left invariant complex product structure
on a Lie group $G$. If $\Gamma\subset G$ is a discrete
subgroup, and $\Gamma$ acts on $G$ by left translations,
then, obviously, (P,J) descends to a complex product structure
on the factor $G/\Gamma$.

Further by a famous theorem of A. Borel \cite{Bore1},
each connected semisimple Lie group $G$ admits uniform
discrete subgroups. The Killing metric (of any signature)
is Einstein (see e.g. \cite{Bes1}). Thus by Theorem \ref{kill1}
the interesting class of examples provided by the next
Corollary is not void.

\begin{coro}\label{cps18}
Let $(P, J)$ be an invariant complex product structure
compatible with the Killing form on $SU(m,m-1)$
(respectively $(Sl(2m-1,\CC))^\RR$ or any simple Lie group).
The factors $SU(m,m-1)/\Gamma$ by any cocompact discrete
subgroup $\Gamma$ are compact complex product manifolds
with compatible non-flat neutral Einstein metric.
\end{coro}

We should remark that non-flat compact complex product manifolds
(See e.g. \cite{IZ1} and references there) have been known for
some time.

 The complex product manifolds introduced in the present paper,
are treated in the context of para quaternionic differential
geometry in \cite{ITZ2}.

\section{preliminaries}\label{prel}

The definitions and standard notation in Lie group and
Lie algebra theory used in this paper can be found e.g.
in \cite{Hel1}. For completenes we review briefly the
definitions and standard facts arround the notion of
complex product structure. For a comprehensive
introduction to the subject we recommend \cite{AnSa1}.

\begin{df}
Let M be a manifold and let $TM$ be the tangent bundle of M.
An {\bf almost product structure} on M is a fibrewise linear
involution, $P:TM \lo TM$ (of constant trace).

An integrable almost product structure on $M$ will be called
{\bf product structure}.

A {\bf complex product structure} on $M$ is a couple $(P,J)$,
where $P$ is a product structure, $J$ is a complex structure,
and we have
\begin{gather}\label{antc}
PJ = -JP.
\end{gather}
\end{df}
Recall that integrability of $P$, respectively $J$, by definition,
means that
\begin{gather}\label{teo4}
\begin{split}
N_P(X,Y) \doteq [P(X),P(Y)] + [X,Y] - P([P(X),Y]) - P([X,P(Y)])= 0;\\
N_J(X,Y) \doteq [J(X),J(Y)] - [X,Y] - J([J(X),Y]) - J([X,J(Y)])=0,
\end{split}
\end{gather}
for any choice of vector fields $X,Y$ on $M$.

It is easy to see that if $P$ is the product structure in a complex
product structure $(P, J)$, then the condition (\ref{antc}) implies
$trP = 0$. This means that if we decompose
the tangent space $TM$ at a point $m\in M$, into eigenspaces
$$
P^+ = P_m^+  \doteq \{X\in T_m:PX = X\},\quad
P^- = P_m^-  \doteq \{X\in T_m:PX = -X\},
$$
then $dimP^+ = dimP^- $ at each point $m\in M$.

\begin{zab}\label{rema4}
Product structures $P$ with zero trace were called {\bf para-complex}
structures by P. Liebermann in the early fifties. They have been
studied by many authors. If $(P, J)$ is a complex product structure
on $M$ and we denote $Q = JP$, then $(Q,J)$ is another complex
product structure on $M$. It is known (in the general case), that the
integrability of $Q$ is a consequence of $(P,J)$ being a complex
product structure (\cite{IZ1},\cite{BVuk0} ). Obviously the situation
is symmetric, we have three anticommuting operators $P, Q, J$ with
\begin{gather}\label{parq}
\begin{split}
P^2 = Q^2 = -J^2 = 1, \\
JP = -PJ = Q,\quad QP = -PQ = J,\quad QJ = -JQ = P.
\end{split}
\end{gather}
The Lie algebra generated by $P, Q, J$ is isomorphic to
$sl(2,\RR)\cong su(1,1)$. The examples constructed in
Section \ref{stru} are such, that the tangent spaces decompose
into real four\footnote{The standard real two dimmensional
representation of $sl(2,\RR)$ also appears in some
examples of CP manifolds, we do not treat this case here.}
dimmensional irreducible representations of $sl(2,\RR)$, all of them
coinciding with the standard representation of $su(1,1)$ in $\CC^2$.
In the literature such complex product manifolds are sometimes called
{\bf hyperparacomplex}\footnote{In algebra and the theory of
arithmetic subgroups, the algebra over $\RR$ generated by $1, J, P, Q$
with (\ref{parq}) satisfied is known as the {\bf indefinite
quaternions}.}. There is a natural class of
pseudoriemannian metrics on such a manifold, which we proceed to define.
\end{zab}

\begin{df}\label{me7}
We shall say that a pseudoriemannian metric $g$ and a
complex product structure $(P, J)$ on a manifold $M$ are
{\bf compatible}\footnote{It is also sensible to
consider metrics such that $g(P(X),P(Y)) = g(X,Y)$,
instead of the condition (\ref{me1}). We do not treat
this case here.} if
\begin{gather}\label{me1}
g(J(X),J(Y)) = g(X,Y),\quad g(P(X),P(Y)) = -g(X,Y),
\end{gather}
for arbitrary tangent vectors $X,Y$ at a point of $M$.
\end{df}
An easy check shows, that a compatible metric has neutral
signature, and that $P^+, P^-$ have to be isotropic subspaces.
A study of such metrics and related connections can be found in
\cite{ITZ2}.

Let $G$ be a simply connected real Lie group of finite dimension,
and let $\mk{g}$ be the Lie algebra of $G$, which we identify
with the space of left invariant vector fields on $G$, thus
trivializing the tangent bundle of $G$.
It is obvious, that a couple of linear operators
$P, J:\mk{g}\lo \mk{g}$, such that
\begin{gather}\label{pcs21}
P^2 = 1,\quad J^2 = -1,\quad PJ = - JP,
\end{gather}
defines a complex product structure on $G$ if and only if
the Nijenhuis conditions (\ref{teo4}) hold for each
$X, Y \in \mk{g}$.

All left translations with elements of the group $G$ acting on
itself are in this case holomorphic transformations of the manifold
$G$ (with respect to the complex structure $J$), and also preserve
the product structure $P$. We say in this case that the complex
product structure $(P, J)$ is {\bf (left) invariant} or more loosely
that $G$ is a {\bf homogeneous complex product}
space\footnote{Obviously the notion of homogeneous $(J, P)$ space
merits a wider definition, we shall not need it.}.

It is well known that right translations on a Lie group $G$
are also holomorphic with respect to a left invariant complex
structure $J$ if and only if
\begin{gather}\label{j1}
[J(X),Y] = J([X,Y]),\quad X,Y \in \mk{g},
\end{gather}
in this case $G$ is a complex Lie group w.r. to $J$. It is also
known, (Proposition 2.6 in \cite{AnSa1}) that only Abelian
Lie algebras admit complex product structures where $J$ satisfies
(\ref{j1}). So we cannot hope to fall in this case with our
simple groups.

In order to make our examlples more transparent, we shall remind
here the alternative definition of integrability. Denote
\begin{gather*}
P^+ \doteq \{X\in \mk{g}: P(X) = X\},
\quad P^- \doteq \{X\in \mk{g}: P(X) = -X\}\\
J^+ = \mk{g}^{1,0} \doteq \{X\in \mk{g}^\CC: J(X) = iX\},\quad
J^- = \mk{g}^{0,1} \doteq \{X\in \mk{g}^\CC: J(X) = -iX\}.
\end{gather*}
Integrability of $(P, J)$ is equivalent to the condition that
$P^+, P^-\quad (J^+, J^-)$ are Lie subalgebras of $\mk{g}\quad
(\mk{g}^\CC)$. Anticommutation of $P, J$
means $J(P^+) = P^-$ (see \cite{AnSa1}).

\section{The complex product structures}\label{stru}
\begin{proof} of Theorem \ref{cps1}

As usual, we denote by $E_j^k\in gl(n)$ the matrix with entry 1 at
the intersection of the j-th row and the k-th column and 0
elsewhere.
     We shall now define two linear operators (P, J) for a
real vector space $\mk{g}$ of dimension $4m^2-4m$ with a base
consisting of elements
\begin{gather*}
U^j, V^j, S^j, T^j, U_j^k, V_j^k, S_j^k, T_j^k,
\end{gather*}
where the range of indices is
\begin{gather}\label{ind}
j = 1,\dots,m -1 ,\qquad j < k < 2m-j.
\end{gather}
 The operators (P, J) are defined for each choice
 of $j, k$ by the conditions:
\begin{gather}\label{mell1}
 P^2 = 1,\qquad J^2 = -1,
\end{gather}
and an explicit formula for each j, k:
\begin{align}\label{sell3}
J(U^j) \doteq V^j;\quad &J(S^j) \doteq T^j;\quad
&P(U^j) \doteq  T^j,\quad &P(V^j) \doteq S^j;\\
\notag
J(U_j^k) \doteq V_j^k;\quad &J(S_j^k) \doteq T_j^k;\quad
&P(U_j^k) \doteq T_j^k; \quad &P(V_j^k) \doteq S_j^k.
\end{align}

We introduce respective bases over $\RR$ for $sl(2m-1,\RR)$ and
$su(m,m-1)$ as real subalgebras of $gl(2m-1,\CC)$.

*We claim that the operators $(P, J)$ defined by formulae
(\ref{mell1}), (\ref{sell3}) and applied to the following
base of $sl(2m-1,\RR)$  define a complex product structure:

\begin{align}\label{sell4}
U^j \doteq &E^j_j + E^{2m-j}_{2m-j} - 2 E^m_m; \quad
&T^j \doteq  E^j_j - E^{2m-j}_{2m-j};\\
\notag
V^j \doteq &E_j^{2m-j} - E_{2m-j}^j;\quad
&S^j \doteq E_j^{2m-j} + E_{2m-j}^j;\\
\notag
U_j^k \doteq  &E_j^k - E_k^j;\quad
&V_j^k \doteq E_k^{2m - j} - E^k_{2m-j};\\
\notag
S_j^k \doteq &E^{2m - j}_k + E^k_{2m-j}; \quad,
&T_j^k \doteq E_j^k + E_k^j.
\end{align}

**We claim that the operators $(P, J)$ defined by formulae
(\ref{mell1}), (\ref{sell3}) and applied to the following
base of $su(m,m-1)$ define a complex product structure:

\begin{align}\label{sell5}
U^j \doteq & i(E^j_j + E^{2m-j}_{2m-j} - 2 E^m_m);
\quad & V^j \doteq i(E^j_j - E^{2m-j}_{2m-j});\\  \notag
T^j \doteq &i(E_j^{2m-j} - E_{2m-j}^j);\quad
&S^j \doteq E_j^{2m-j} + E_{2m-j}^j.
\end{align}
\begin{equation}\label{sell15}
\begin{split}
U_j^k \doteq &
\begin{cases}
E_j^k - E_k^j& \text{ if }\quad  j< k \leq m; \\
E_k^{2m - j} - E_{2m-j}^k & \text{ if }\quad m < k < 2m-j.
\end{cases}\\
V_j^k \doteq &
\begin{cases}
i(E_j^k + E_k^j)& \text{ if }\quad  j< k \leq m;\\
i(E_k^{2m - j} + E_{2m-j}^k) & \text{ if }\quad m < k < 2m-j.
\end{cases}\\
S_j^k \doteq &
\begin{cases}
E_k^{2m - j}+ E_{2m-j}^k& \text{ if }\quad  j < k \leq m; \\
- E_j^k - E_k^j & \text{ if }\quad m < k < 2m-j.
\end{cases}\\
T_j^k \doteq &
\begin{cases}
i(E_k^{2m - j}- E_{2m-j}^k)& \text{ if }\quad  j < k \leq m; \\
 i(E_k^j - E_j^k)& \text{ if }\quad m < k < 2m-j.
\end{cases}
\end{split}
\end{equation}

The proof of claims * and ** is a straightforward check,
using the integrability conditions (\ref{teo4}). One has to check
only the case $j = 1,\quad 2 \leq k \leq 2m-2$.
Indeed,
let $\mk{g}$ be either $su(m,m-1)$ or $sl(2m-1,\RR)$.
The operators $P,J$, preserve the subspace
\begin{gather}\label{jo2}
span\{U^1, V^1, S^1, T^1, U_1^k, V_1^k, S_1^k, T_1^k:
1 <k < 2m-1\}.
\end{gather}
Thus if we remove the external two rows and columns (\ref{jo2}) in
all matrices of $su(m,m-1)\quad(sl(2m-1,\RR))$, then we get the Lie
subalgebras
$$
su(m-1,m-2) \subset su(m,m-1),\quad (sl(2m-3,\RR)\subset sl(2m-1,\RR)).
$$
By their definition, the operators $P, J$ preserve these
subalgebras and restricting to them we get the induction step.
\end{proof}

\begin{zab}\label{rema6}
In the last paragraph of the proof above, we have shown, that
our complex product structures allow embeddings
\begin{gather*}
SU(2,1) \subset SU(3,2)\subset\dots\subset SU(m,m-1);\\
SL(3) \subset SL(5)\subset\dots\subset SL(2m-1)
\end{gather*}
as complex product submanifolds. Similar induction was
used by Joyce \cite{Joyce2} to get
hypercomplex structures on compact reductive groups.
\end{zab}

Using the symmetry of $(P, J)$ defined above, we can get some
further examples

\begin{zab}\label{rema1}
Let $\mk{g}$ be either $su(m,m-1)$ or $sl(2m-1,\RR)$.
The operators $(P, J)$ of the proof of Theorem \ref{cps1}
preserve the subspace
\begin{gather*}
span\{E_m^j,E_k^m: 1 \leq j,k \leq 2m-1\}\cap\mk{g}
\end{gather*}

Thus, if we remove the middle row and column in all
matrices of $su(m,m-1)$ (resp. $sl(2m-1,\RR)$), then we get the
subalgebras
$$
u(m-1,m-1) \subset su(m,m-1) \quad (\text{ respectively }
gl(2(m-1),\RR)\subset sl(2m-1,\RR))
$$
The action of (P, J) preserves these subalgebras also, so for each
$m\geq 1$ we have found a complex product structure on $U(m,m)$
( $GL(2m,\RR)$)  and embedded it as a complex
product submanifold  in $SU(m+1,m)$ (respectively $SL(2m+1,\RR)$).
\end{zab}

We shall now introduce parameters in the definition of $J, P$.
Let $\mk{g}$ be one of the algebras $su(m,m-1)$ or $sl(2m-1,\RR)$.
Using the notation of the proof of Theorem \ref{cps1} we define
an Abelian subalgebra
\begin{gather}\label{ki}
\mk{z} \doteq \RR{U^1}\oplus\dots\oplus\RR{U^{m-1}} \subset \mk{g}.
\end{gather}
We have

\begin{prop}\label{ki2}
Let $\mk{g}$ be any of the Lie algebras $su(m,m-1), sl(2m-1,\RR)$. Let
$Z^1,\dots,Z^{m-1}$ be any base for the subalgebra
$\mk{z}$ defined in (\ref{ki}). Let
\begin{gather*}
V^j, S^j, T^j, U_j^k, V_j^k, S_j^k, T_j^k,
\end{gather*}
be as in the proof of Theorem \ref{cps1}.
Then the operators $J, P:\mk{g}\lo\mk{g}$ defined by
\begin{align}\label{sell77}
\notag  P^2 = -J^2 = 1\\
J(Z^j) \doteq V^j;\quad &J(S^j) \doteq T^j;\quad
&P(Z^j) \doteq  T^j,\quad &P(V^j) \doteq S^j;\\
\notag
J(U_j^k) \doteq V_j^k;\quad &J(S_j^k) \doteq T_j^k;\quad
&P(U_j^k) \doteq T_j^k; \quad &P(V_j^k) \doteq S_j^k.
\end{align}
give a complex product structure on $\mk{g}$.
\end{prop}
\bpr
We have to prove integrability. Again a straightforward check of the
Nijenhuis conditions gives the result. It helps to check first that
$$
[Z,J(X)] = J([Z,X]);\quad [Z,P(X)] = P([Z,X]),\qquad Z\in\mk{z},
X\in \mk{g}.
$$
\epr

Now we turn to metric properties of the above 
complex product structures.

\begin{proof} of Theorem \ref{kill1}

The scalar product
\begin{gather}\label{scal}
\scal{X,Y} \doteq \frac{1}{2}tr(XY),\quad X,Y \in su(m,m-1),
\end{gather}
is proportional to the Killing form.

First we treat $su(m,m-1)$. The product (\ref{scal}) is negative
definite on the subalgebra $\mk{z}$ (defined in (\ref{ki})). Let
$Z^1,\dots,Z^{m-1}$ be any orthonormal base of $\mk{z}$ with respect
to $\scal{.,.}$. We take $ V^j, S^j, T^j, U_j^k, V_j^k, S_j^k, 
T_j^k$, as defined in formulae (\ref{sell5}), (\ref{sell15}). 
Then the base
\begin{gather}\label{base}
Z^j, V^j, S^j, T^j, U_j^k, V_j^k, S_j^k, T_j^k,
\end{gather}
is orthogonal with respect to $\scal{.,.}$ and
\begin{gather*}
\norm{Z^j}^2 = \norm{V^j}^2 = \norm{U_j^k}^2 = \norm{V_j^k}^2 = -1;\\
\norm{S^j}^2 = \norm{T^j}^2 = \norm{S_j^k}^2 = \norm{T_j^k}^2 =  1.
\end{gather*}
It is then obvious from the definition of $P, J$ (in formula
(\ref{sell77})), that $\scal{.,.}$ is compatible, thus $su(m,m-1)$
is settled.

We present $sl(2m-1,\CC) = su(m,m-1) \oplus isu(m,m-1)$. Then
we keep the action of $(P, J)$ on  $su(m,m-1)$ as in Proposition 
\ref{ki2} with respect to the base (\ref{base}),
and extend them to $sl(2m-1,\CC)$ by
\begin{align}\label{rr}
J(iZ^j) \doteq iV^j;\quad &J(iS^j) \doteq iT^j;\quad
&P(iZ^j) \doteq  iT^j,\quad &P(iV^j) \doteq iS^j;\\
\notag
J(iU_j^k) \doteq iV_j^k;\quad &J(iS_j^k) \doteq iT_j^k;\quad
&P(iU_j^k) \doteq iT_j^k; \quad &P(iV_j^k) \doteq iS_j^k.
\end{align}
The scalar product $\scal{X,Y}$ defined in (\ref{scal})
(now we apply it to $X,Y\in sl(2m-1,\CC)$) is
proportional to the Killing form on $sl(2m-1,\CC)$. Then
\begin{gather*}
K(X,Y) \doteq \Re\scal{X,Y}
\end{gather*}
is proportional to the Killing form on $sl(2m-1,\CC)^\RR$
(see e.g. \cite{Hel1}, Ch. III, Lemma 6.1). Thus the set of
\begin{gather*}
Z^j, V^j, S^j, T^j, U_j^k, V_j^k, S_j^k, T_j^k,\\
iZ^j, iV^j, iS^j, iT^j, iU_j^k, iV_j^k, iS_j^k, iT_j^k.
\end{gather*}
with any indices satisfying (\ref{ind}), is an orthogonal base of
$sl(2m-1,\CC)^\RR$ and
\begin{gather*}
K(iZ,iZ) = K(iU,iU) = K(iV,iV) = K(S,S) = K(T,T) = 1,\\
K(Z,Z) = K(U,U) = K(V,V) = K(iS,iS) = K(iT,iT) = -1.
\end{gather*}
Whence the extention (\ref{rr}) of $(P, J)$ to $sl(2m-1,\CC)^\RR$ is
compatible with the Killing metric.
\end{proof}

\begin{ack}
The research is partially supported by Contract MM 809/1998 with
the Ministry of Science and Education of Bulgaria, Contracts
586/2002 and 35/2003 with the Sofia University "St. Kl. Ohridski". S.I. is
a  member of the EDGE, Research Training Network
HPRN-CT-2000-00101, supported by the European Human Potential
Programme.
\end{ack}


\begin{thebibliography}{99}

\bibitem{AnSa1} Andrada A., Salamon S., Complex product
structures on Lie algebras, Preprint DG/0305102.

\bibitem{Bes1} Besse A., Einstein manifolds, Springer, 1987

\bibitem{BVuk0}  Blazic N., Vukmirovic S., Four-dimensional Lie
algebras with parahypercomplex structure, math. DG/0310180.

\bibitem{Bore1} Borel A., Compact Clifford-Klein forms of
symmetric spaces, Topology, 2(1963), 111-122

\bibitem{Hel1} Helgason S., Differential Geometry, Lie Groups, and
Symmetric Spaces, Academic Press, 1978

\bibitem{ITZ2} Ivanov S., Tsanov V., Zamkovoy S.,
Hyper-Parahermitean manifolds with torsion, math.DG/0405585

\bibitem{IZ1} Ivanov S., Zamkovoy S., Para-Hermitian and
Para-Quaternionic manifolds, math.DG/0310415.

\bibitem{Joyce2} Joyce D., Compact hypercomplex and quaternionic
manifolds, J. Diff. Geom., 35(1992), 743 - 761
\end{thebibliography}
\end{document}